\documentclass[a4paper]{article}
\pdfoutput=1
\usepackage{amsmath}
\usepackage{amssymb}
\usepackage{layout}
\usepackage{ifthen}
\usepackage{url}
\usepackage[pdftex,colorlinks=true,linkcolor=blue]{hyperref}
\usepackage[retainorgcmds]{IEEEtrantools}
\author{Roger Sewell}
\usepackage{graphicx}
\usepackage{xr}
\usepackage{float}

\title{Assessment of the quality of a prediction\footnote{Fifth
  version deposited in arxiv; RFS version 1.19.1.2
; corrected reference numbering and missing details for reference 13.}}

\author{Roger Sewell\footnote{Trinity College, Cambridge; but at the
  time this work was done employed full time by Cambridge Consultants.}\\
\href{mailto:roger.sewell@cantab.net}{\scriptsize{roger.sewell@cantab.net}}
}

\DeclareMathOperator{\E}{E}
\DeclareMathOperator{\trace}{trace}
\DeclareMathOperator{\mean}{mean}
\DeclareMathOperator{\cov}{cov}

\begin{document}
\maketitle

Acknowledgements: To Elisabeth Crowe, David Rimmer, Edwin Carter, and
Sam Pumphrey for helpful comments.

\begin{center}
\textbf{Abstract}
\end{center}
Shannon defined the mutual information between two variables. We
illustrate why the true mutual information between a variable and the
predictions made by a prediction algorithm is not a suitable measure
of prediction quality, but the apparent Shannon mutual information
(ASI) is; indeed it is the unique prediction quality measure with
either of two very different lists of desirable properties, as
previously shown by de Finetti and other authors. However, estimating
the uncertainty of the ASI is a difficult problem, because of long and
non-symmetric heavy tails to the distribution of the individual values
of $j(x,y) = \log \frac{Q_y(x)}{P(x)}$. We propose a Bayesian modelling
method for the distribution of $j(x,y)$, from the posterior
distribution of which the uncertainty in the ASI can be inferred. This
method is based on Dirichlet-based mixtures of skew-Student
distributions. We illustrate its use on data from a Bayesian model for
prediction of the recurrence time of prostate cancer. We believe that
this approach is generally appropriate for most problems, where it is
infeasible to derive the explicit distribution of the samples of
$j(x,y)$, though the precise modelling parameters may need adjustment
to suit particular cases.

\tableofcontents

\renewcommand{\arraystretch}{1.5}

\section{Introduction}

In many machine-learning situations a prediction or classification
algorithm is trained on data for which the right answer is known -- so
called ``supervised'' learning on training data – and then required to
make predictions on data to which it has not previously been
exposed. There are numerous methods for assessing the quality of the
resulting predictions; many depend on a particular measure of cost of
error. For example, where the prediction provides a point value rather
than a probability distribution, the expectation of squared difference
between the true value and the predicted value is one such measure.

However, there are many situations where a variety of cost functions
may be relevant to prediction, depending on who is using the
prediction\footnote{For an example to keep in mind while reading this
paper, consider a patient who has been diagnosed with cancer, and is
interested in when he will die, in order that he may thriftily take
out life insurance only for the periods in which the risk is high. The
insurance company, his doctor, the patient himself, and his relatives
may all have very different views on what cost function is
appropriate, but all will benefit from having an accurate probability
distribution on time of death from a prediction algorithm.}. We adopt
the Bayesian paradigm, so that a high quality but imperfect prediction
will provide not a point value of the variable being predicted, but a
probability distribution on it\cite{Lindley,
  DJCMbook}. Hilden\cite{Hilden} surveyed many ways of measuring the
quality of such predictions.

It is our belief that one method in particular provides the most
appropriate means of assessing the performance of the upper end of the
quality range for prediction algorithms independent of any specific
cost function; that measure is the Apparent Shannon Information (ASI)
in the predictions about the true values. This paper recalls that
method and its desirable properties. We claim no originality in
inventing this method (believing Claude Shannon\cite{Shannon} to be
its conceptual father in the 1940s), nor in realising the nature of
its desirable properties\cite{Winkler, Shapiro, DJCMpersonal} (with
minor differences). We strongly believe that it is used far less often
than it should be, and therefore wish to remind those constructing and
assessing prediction algorithms of its benefits.

However, there are difficulties assessing the value of this measure
using only finite amounts of data. We also offer in this paper a
Bayesian approach to evaluating the ASI, based on skew-Student mixture
models, and illustrate the measurement with results from work on
prediction of the recurrence of prostate cancer, also using
skew-Student mixture models.

We also do not claim that this method is suitable for distinguishing
the difference in performance of two imperfect prediction algorithms
that produce only point-value outputs rather than probability
distributions over the possible outputs.

Further, we do not consider in this paper the cost of execution of any
algorithm, merely the quality of the predictions it makes. Clearly in
any real-life choice of algorithm, such issues enter the decision
along with the quality of the predictions made, and if formal
trade-offs between error costs and implementation costs are needed it
will be necessary to specialise further to the particular cost
function that is relevant.

The rationale for using the ASI and definitions are in section
\ref{rationale}. Section \ref{properties} recalls the properties of
the ASI. Section \ref{stability} illustrates the need for a means of
estimating the ASI beyond simply averaging log probability ratios over
a finite unseen dataset, and provides that means, and section
\ref{conclusions} presents our conclusions.

\section{Definitions and rationale}
\label{rationale}

\subsection{Notation}

We suppose we wish to predict the value of a random variable $x$
taking values in some measurable space $X$, and to do so on the basis
of some data $y$ taking values in a measurable space $Y$. We suppose
also that we have some prior information on the likely values of $x$
in the absence of knowledge of $y$. Given a possible prediction
algorithm, we suppose it in general to output a probability
distribution $Q_y$ on $x$ which depends on $y$; in general this
distribution may not be the Bayesian posterior distribution $P(x|y)$,
but as a probability distribution it must none the less have integral
of $1.0$. We will understand this to cover also the situation that the
algorithm outputs a point value of $x$, when for any measurable subset
$X_1 \subseteq X$, $Q_y(X_1)$ will be $1$ if $X_1$ contains that point
value of $x$ and $0$ otherwise.

However, for simplicity and ease of general understanding we will
write formulae as if all variables are continuous with an appropriate
probability density, so instead of $\int{f(x)\,dP(x)}$ we will write
$\int{f(x)P(x)\,dx}$, supposing also that $Q_y(x)$ denotes the output
probability density on $x$, and leave the reader to make the necessary
notational adjustments for (partially) discrete random variables with
non-zero probabilities on single-point sets.

\subsection{Reasoning leading to choice of Apparent Shannon Information,
definitions, and basic properties}

Shannon\cite{Shannon} defined the mutual information $I(x;y)$ between
$y$ and $x$ by $$I(x;y)=\int_{X\times
  Y}{P(x,y)\log\left(\frac{P(x|y)}{P(x)}\right)\,d(x,y)};$$ to
distinguish it from the Apparent Shannon Information introduced later,
we will refer to $I(x;y)$ as the true Shannon mutual information in
$y$ about $x$. If the logarithm is taken to the base $2$ the resulting
value is reported in bits; if, as we shall, the logarithm is taken to
the base $e$ the resulting value is reported in nats or nepers (in
this context the two have the same meaning).

We first need to see why $I(x;Q_y)$ is not a suitable measure of a
prediction algorithm’s quality. Suppose $x$ is a variable taking only
two possible values (e.g. whether a coin will come down heads or tails
when next thrown), and that in the absence of any data each of the two
values is equally probable. Then any algorithm which tells us the
right answer with probability 1 will give us 1 bit, or about 0.69
nepers, of true mutual information in its output about $x$. However, so
will any algorithm which with probability 1 tells us the wrong answer;
this is because an answer which is known to be wrong carries in this
situation as much relevant information as one which is known to be
right, as can be seen by applying an inverter to the algorithm’s
output. It is immediately clear to any betting man that an oracle
which can tell which side of a coin will not be uppermost is in some
senses just as good as one which can tell which side will be visible.

However, we are interested in assessing the quality of the original
algorithm, without any correction mechanism added to its output. We
therefore define the Apparent Shannon Information by $$J(x;Q_y) =
\int_{X\times
  Y}{P(x,y)\log\left(\frac{Q_y(x)}{P(x)}\right)\,d(x,y)},$$ where the
substitution of $Q_y(x)$ for $P(x|y)$ replaces what can be deduced
from the output of the algorithm with what the output of the algorithm
is actually claiming. Further, we consider also a generalisation of
this, where before developing the present algorithm $A_1$ we have a
perhaps inferior algorithm $A_0$, and we want to measure the quality
of $A_1$ relative to that of $A_0$. We generalise and
define $$J(x;Q_y^{A_1};Q_y^{A_0}) = \int_{X\times
  Y}{P(x,y)\log\left(\frac{Q_y^{A_1}(x)}{Q_y^{A_0}(x)}\right)\,d(x,y)},$$
the apparent Shannon information from $A_1$ relative to that from
$A_0$, but where there is an understood prior distribution on $x$ we
will continue to write $J(x;Q_y)$. Note that the ASI corresponds to
the difference of expectations of log probability assigned to the true
values by the algorithm and by the prior (or the antecedent
algorithm). As such the ASI differs only in setting of the zero point
from the expected log probability measures of \cite{Winkler} and
\cite{Shapiro} (which also differ in zero point from each other).

Now, while it is always the case that both $I(x;y)\geq 0$ and
$I(x;Q_y)\geq 0$, it is not always the case that $J(x;Q_y)\geq 0$;
where $J(x;Q_y) < 0$, the output distribution $Q_y$ is misleading in
comparison with the prior (or antecedent algorithm's)
distribution. The example above where an algorithm reliably predicts
the wrong outcome of a coin-toss is typical of a misleading output; it
scores $J(x;Q_y) = -\infty$. Equally, any other algorithm in any
situation which outputs a point value (or equivalently all the
probability on a single point value) will score $J(x;Q_y) = -\infty$ if
the probability of that point value being wrong is non-zero.

\subsection{Conditions for measurement}
We wish to emphasise particularly that $J(x;Q_y)$ is to be calculated
with the distribution $P(x,y)$ reflecting the distribution of new
incoming problems, not some version of this that has a distribution
biased in favour of seeing cases on which an algorithm has been
``trained''. In some situations there will be large differences
between values calculated on ``unseen'' data (as they should be) and
values calculated based on data seen in a training set. Therefore when
calculating $J(x;Q_y)$ data should only be used that has not been seen
during training of the algorithm.

In particular Bootstrapping (\cite{Harrell} page 372) is not an
appropriate way of determining $J(x;Q_y)$, as it does not correctly
calculate the differences between cases seen and unseen during
training (see appendix to \cite{prostate}).  

We will discuss methods of estimating the true value of $J(x;Q_y)$ in
section \ref{stability} below.

\section{Properties}
\label{properties}

We now briefly recall the elementary properties of this prediction
quality measure.  These follow easily from the Kullback-Leibler
divergence theorem, also known as Gibbs’ inequality (\cite{DJCMbook}
pp 34 ff), and from the Data Processing Theorem\cite{Cover,
  ShannonWeaver}. (For this we assume that the variable being
predicted is a real scalar, though this is not an essential
limitation.)

\subsection{List of basic properties}
\label{proplist}

\begin{enumerate}

\item \label{Bayesmaxes} The true Bayesian posterior distribution
  $Q_y^{\text{Bayes}}(x) = P(x|y)$ maximises $J(x,Q_y)$ over all
  probability distributions $Q_y$ , at the value $J(x;
  Q_y^{\text{Bayes}}) = I(x;y)$.

\item We have the following relationships between the ASI and true
  Shannon mutual information: $J(x;Q_y) \leq I(x;Q_y) \leq I(x;y)$. In
  other words the ASI never exceeds the true Shannon mutual
  information in the prediction, which in turn never exceeds the
  mutual information between data and predicted variable.

\item \label{transformability} If $w=f(x)$ for some bijective
  differentiable function $f$ with differentiable inverse, then
  $J(w;R_y) = J(x;Q_y)$, where $R_y(w) = Q_y(g(w))g'(w)$, where $g'$
  denotes the derivative of $g$, the inverse of $f$. For example, if
  we can predict $x$, a positive real variable, with a particular
  amount of ASI, then we can also predict $\log(x)$ with the same ASI.

\item \label{serialadditivity} If we start with a prior $P(x)$, then
  develop algorithm $A_1$, then develop an improved algorithm $A_2$,
  the ASI of $A_2$ relative to the prior is simply the sum of that of
  $A_1$ relative to the prior and that of $A_2$ relative to $A_1$:
  $J(x; Q_y^{A_2} ; Q_y^{A_1}) + J(x; Q_y^{A_1}) = J(x;Q_y^{A_2}).$

\item \label{eventadditivity} If there are two mutually exclusive and
  exhaustive events $E_0, E_1$ which are not observed, then $J(x;Q_y)$
  will be the appropriate weighted average of $J_{E_0}(x;Q_y)$ (the
  ASI when $E_0$ occurs but is not observed) and $J_{E_1}(x;Q_y)$ (the
  ASI when $E_1$ occurs but is not observed), i.e. $J(x;Q_y) =
  P(E_0)J_{E_0}(x;Q_y) + P(E_1)J_{E_1}(x;Q_y)$.

\item \label{unique1} The ASI is \cite{Winkler, Shuford, deFinetti},
  up to a constant scalar multiple, the unique quality measure that is
  zero on the prior, is proper in the sense that the person making the
  prediction has no incentive to mis-state his beliefs, and depends
  only on $Q_y$.

\item \label{unique2} The ASI is also, up to a constant scalar
  multiple, the unique quality measure that is zero on the prior, and
  has properties \ref{transformability}, \ref{serialadditivity}, and
  \ref{eventadditivity} above (outline proof in appendix section
  \ref{appxproof} below); the scalar multiple can be resolved by using
  also property \ref{Bayesmaxes}.

\end{enumerate}

\subsection{Betting outcome prediction}
\label{betting}

One may also ask what the significance of a positive value of
$J(x;Q_y)$ is. One possible characterisation is via the following
betting scheme.

Each of two individuals 1 and 2 has a prediction algorithm ($A_1$ and
$A_2$ respectively) for predicting $x$ from $y$. Each player places a
sum $M$ of money into a pile in front of him, keeping the two piles
separate. Each time a new value of $y$ is observed, each player $i$
makes his prediction $Q_y^{A_i}$, a distribution on $x$. Each time the
true value $x$ corresponding to a prediction becomes known, the
amount of money currently in each player $i$'s pile is multiplied by
$\frac{Q_y^{A_i}(x)}{P(x)}$. At the end of the game each player
pockets the pile in front of him.

It is then easy to verify that the geometric expected value in each
pile $i$ approaches $\infty$ as the number of predictions checked
approaches $\infty$ if and only if the predictions being made satisfy
$J(x;Q_y^{A_i}) > 0$; the geometric expected ratio of player 1’s pile
to player 2’s pile approaches $\infty$ if and only if $J(x;Q_y^{A_1})
> J(x;Q_y^{A_2})$. In other words betting according to the predictions
works without diversification if and only if the predictions carry
positive apparent Shannon information.

\subsection{Relevance to design of algorithms based on Bayesian
  inference} 

Algorithms based on Bayesian inference can be designed using the
inherent capability of Bayesian techniques to consider a range of
different models for the probability distributions involved, and to
choose optimally and probabilistically between them. So why then might
one want to measure the performance of such an algorithm ?

Consider a situation where pairs $(x,y)$ arrive to be processed, with
$y$ observed but $x$ hidden and to be inferred or predicted. The
population arriving comes from a joint distribution $P(x,y)$; suppose
that the true marginal distribution $P(x)$ and the true likelihood
$P(y|x)$ are known, so that the true Bayesian posterior $P(x|y)$ can
be calculated.

It would then be appropriate that the prediction quality measure is
optimised by setting $Q_y(x) = P(x|y)$. This is indeed the case, as
seen in section \ref{proplist} item \ref{Bayesmaxes} above.

However, in many situations a Bayesian investigator is handicapped by
only knowing his own prior on the unknown and his own approximation to
the likelihood. Either or both of these may have been learned from
data, for example training data. Such learning may well have been done
by comparing the posterior model probabilities for a number of
possible models (distribution families), and the choice between models
may have been done optimally, given the available
information. Nonetheless, the question of which initial candidate set
of models to choose between is one to which Bayesian techniques
provide no definitive answer. Therefore it is still of interest, given
two Bayesian algorithms designed from different selections of possible
models, to ask which is better.

One way of doing this is to combine the sets of training data used in
each development and use the combined set in the standard Bayesian
model choice procedure with both models included on the candidate
list; indeed this is expected to produce a combined algorithm at least
as good as the better algorithm and possibly better than both.
Sometimes, however, access to one of the training sets or model
definitions is limited for e.g. commercial reasons, and in that
situation the measurement of $J(x;Q_y)$ for each algorithm provides an
appropriate answer. One must, however, ensure that data previously
seen by either algorithm is excluded from such a test set.

\section{Stability and estimation}
\label{stability}

Given a prediction algorithm and some unseen data with corresponding
right answers, how do we measure the apparent information content ?

For example, we might develop an algorithm for prediction of time of
recurrence of prostate cancer, based on a training dataset of
information about a series of patients and the times cancer recurred
in those patients. We then might have a disjoint set of patients for
whom we have the same input information $y$ and actual outcome data
$x$, for whom we produce predictions $Q_y$, and we want to measure
$J(x;Q_y)$.

\subsection{A problem: point estimates based on finite datasets do not
  suffice}

Given that $J(x;Q_y) =
\int{P(x,y)\log\left(\frac{Q_y(x)}{P(x)}\right)\,d(x,y)} = \E
\log\left(\frac{Q_y(x)}{P(x)}\right)$, one obvious approach is to
observe previously unseen data $y$ and corresponding correct answers
$x$, apply the $Q_y(x)$ prediction algorithm, and calculate for each
$(x,y)$ pair the value of $j(x,y) =
\log\left(\frac{Q_y(x)}{P(x)}\right)$, then take the average of these
values over the data points observed. (It might even be tempting to
suggest that the uncertainty in the result would be Gaussian with
standard deviation equal to that of $j(x,y)$ divided by the square
root of the number of points observed minus one. This, however, would
be very misleading.)

\begin{figure}
\begin{center}

\includegraphics[scale=0.6]{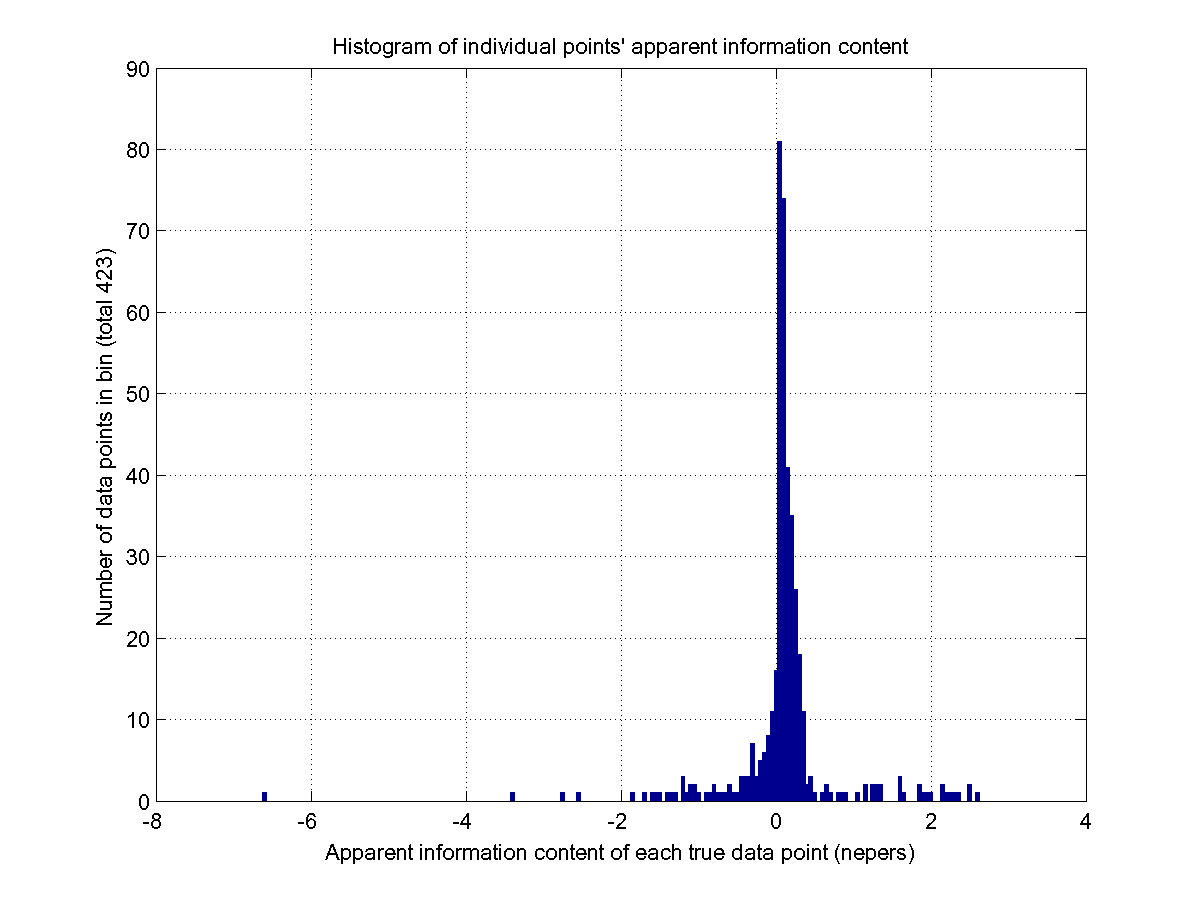}

\end{center}

\caption{Histogram of $j(x,y)$ for data points not seen during
  training, where the prediction algorithm is using a range of
  biomarkers and clinical data to predict time of recurrence of
  prostate cancer following radical prostatectomy.  
\label{fig1}
}

\end{figure}

The essential problem is illustrated by the measurements of $j(x,y)$
taken from just such a Bayesian algorithm trained to predict time of
recurrence of prostate cancer from a set of biomarkers\cite{prostate},
a histogram of which on data points not seen during training is shown
in Figure \ref{fig1}. The average value of $j(x,y)$ over these data
points is $+0.096$ nepers, but if the one point at $-6.58$ nepers were
not included it would be $+0.112$, which given the narrowness of the
main distribution is quite a big difference (the false assumption that
these data were Gaussianly distributed would suggest that the mean is
known to a standard deviation of $0.034$ nepers, making this change
more than $0.5$ standard deviations).

We therefore need to know not only the average of the $j(x,y)$, but
its distribution, in order that we may know just how uncertain our
quality measure is.

\subsection{The suggested approach}
\label{approach}

Rasmussen and others have investigated using Student distributions to
model this data and infer uncertainty\cite{DJCMpersonal, DELVE,
  DELVEtechreport}, encountering difficulty because of the tendency
for there to be long heavy tails that are not symmetric. We propose a
related method, using the skew-Student distribution, which lacks
reflection symmetry about its mean, as follows.

We create a parametric Bayesian model for the distribution of
$j(x,y)$, then draw random samples $(\theta_k)_{k=1,...,K}$ from the
distribution on the parameters given the $j(x,y)$ data. For each
sample $\theta_k$ we can calculate the mean $\mu_k$ of the resulting
sample of the distribution, and hence get a histogram, cumulative
distribution, etc of the desired value $J(x;Q_y)$.

\begin{figure}
\begin{center}

\includegraphics[scale=0.25]{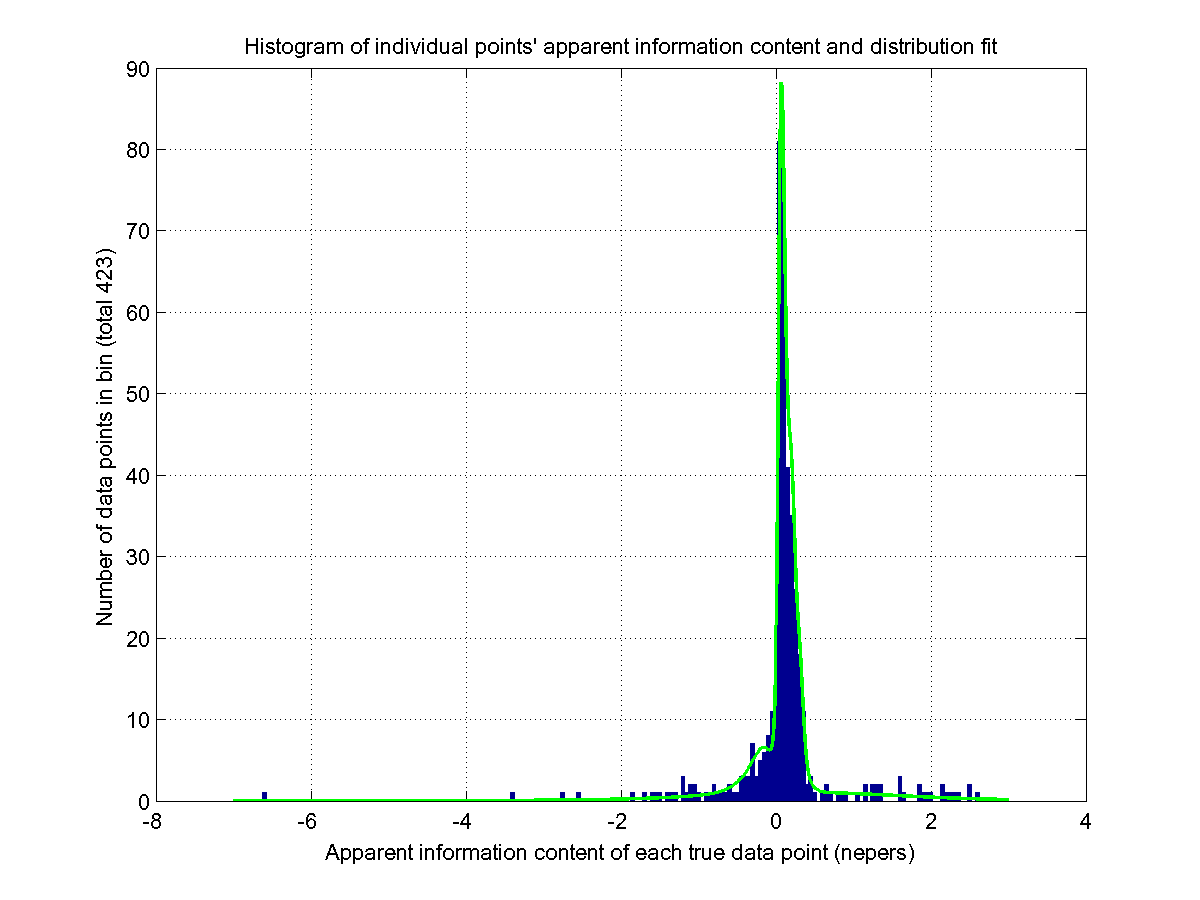}
\includegraphics[scale=0.25]{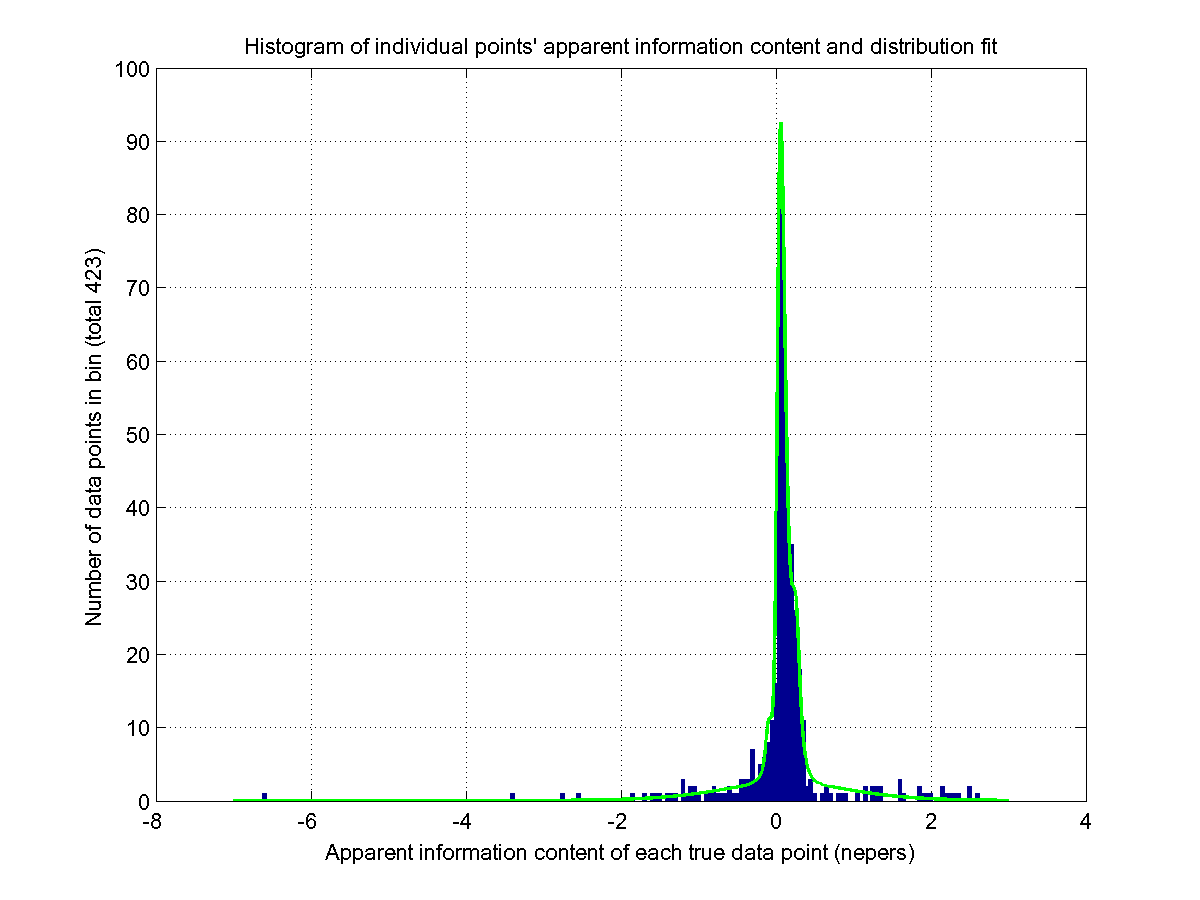}
\includegraphics[scale=0.25]{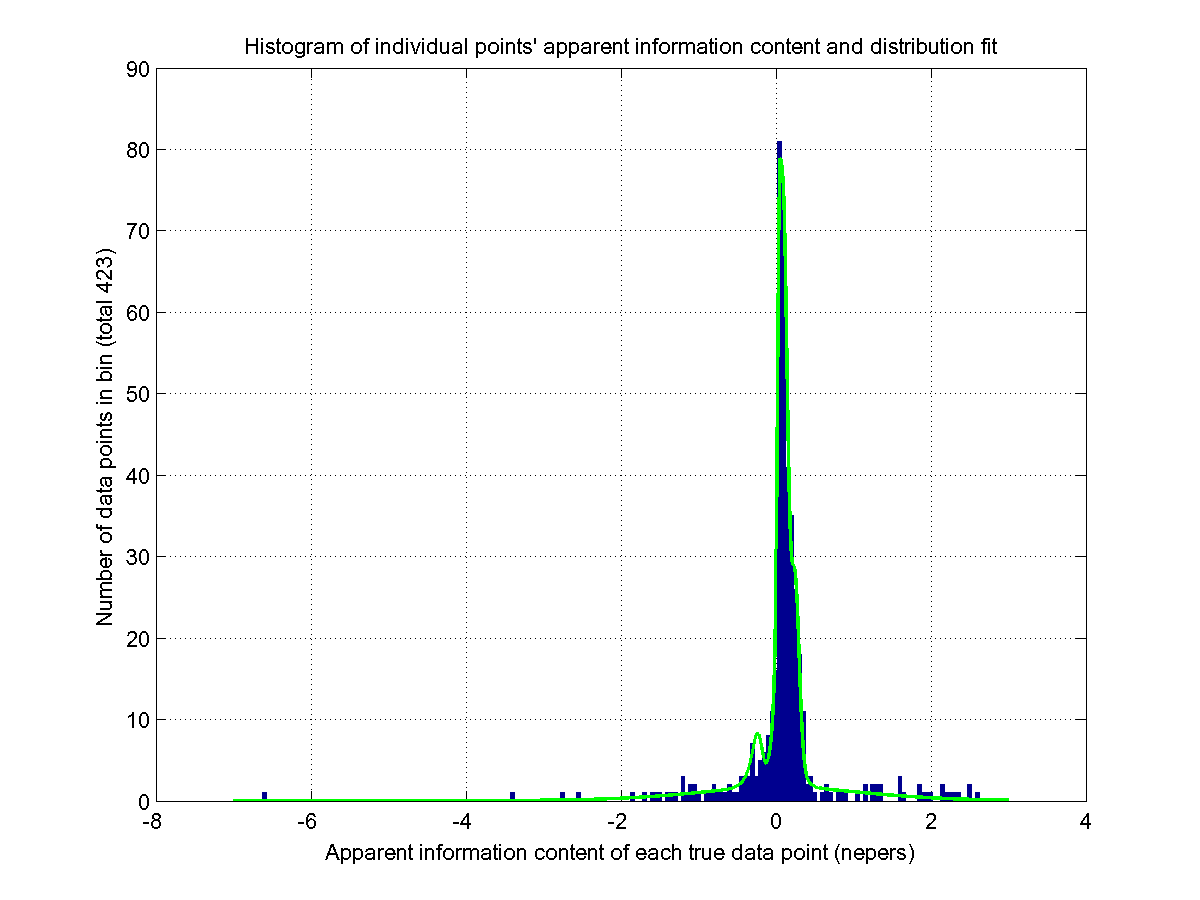}

\end{center}

\caption{Three samples of the distributions specified by the posterior
  distribution of the $\theta_k$ given the $j(x,y)$ data. For each
  (green) distribution a mean value can be calculated, since the
  distribution’s parameters are known. Given a large number of such
  samples of the mean, the posterior distribution of the mean can be
  reconstructed as in Figure \ref{fig3}.
\label{fig2}
}

\end{figure}

\begin{figure}
\begin{center}

\includegraphics[scale=0.6]{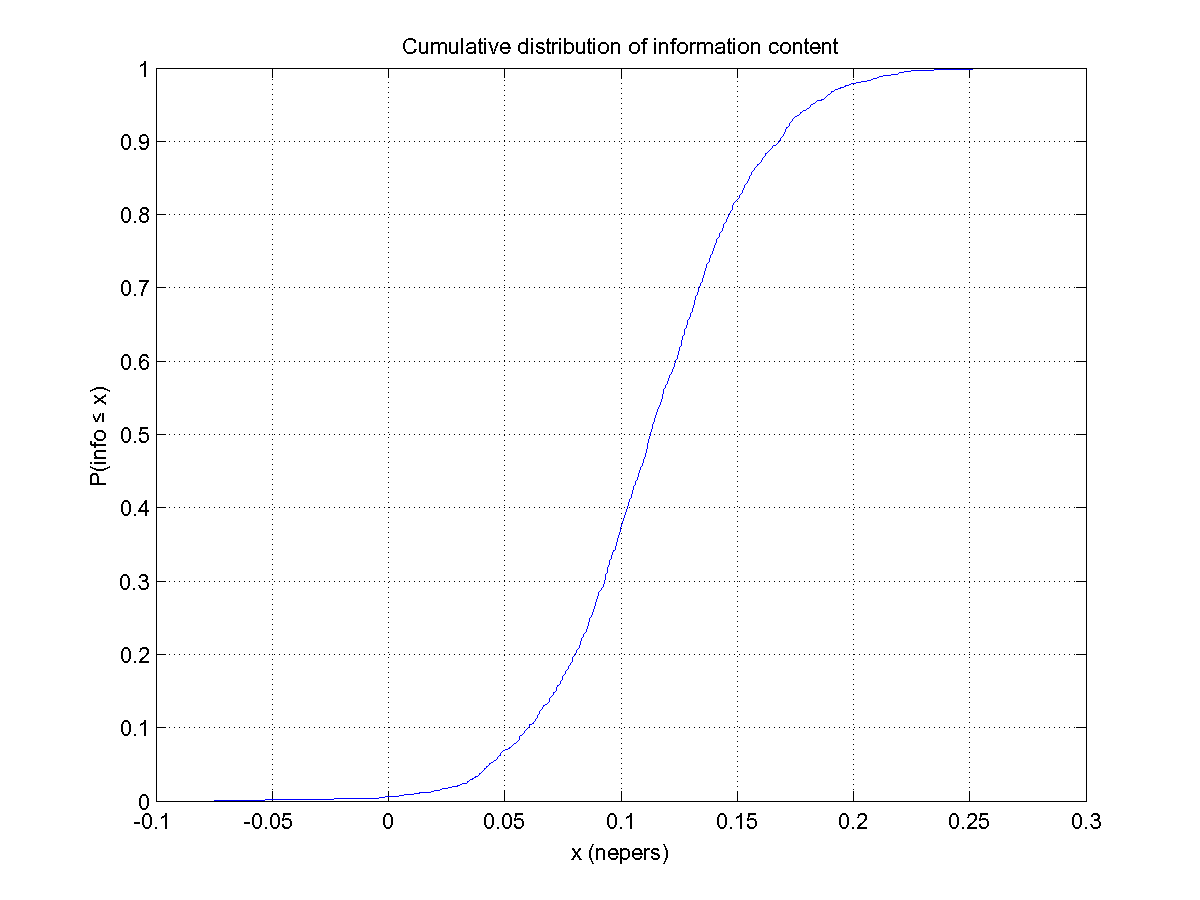}

\end{center}

\caption{The cumulative distribution of $J(x; Q_y)$ based on modelling
  the $j(x,y)$ data calculated from the data points unseen during
  training.
\label{fig3}
}

\end{figure}

We illustrate how this idea works before defining the details of the
modelling. Figure \ref{fig2} shows the actual distributions resulting
from three of the samples of $\theta_k $ when this modelling process
is run on the $j(x,y)$ data of Figure \ref{fig1}, and Figure
\ref{fig3} shows the resulting cumulative posterior distribution of
$J(x;Q_y)$ given the $j(x,y)$ data.

This cumulative distribution and the samples that gave it allow us to
report our uncertainty on $J(x;Q_y)$ for example as mean $0.107$,
median $0.106$, $0.025$ quantile $0.043$, and $0.975$ quantile $0.172$
nepers. The specific value $0.096$ which is the mean of the $j(x,y)$
not surprisingly lies comfortably with these results.

Obviously the model used for this purpose needs to be appropriately
flexible. The method used here is to create a hierarchical Bayesian
model based around a Dirichlet-mixed skew-Student mixture; the
details are in the appendix section \ref{modelling}. The key element
of this is the skew-Student distribution, which is the marginal
distribution of $\frac{\sigma z+\nu}{\sqrt{\alpha}}+\mu$ that results
from $z$ being distributed according to a standard unit Gaussian and
$\alpha$ being Gamma distributed independently of $z$, where $\sigma,
\mu, \nu$ are parameters, along with $m$ and if desired $r$,
respectively the shape and scale parameters of the relevant Gamma
distribution.

Obviously it is important that in the absence of $j(x,y)$ data the
resulting distribution is sufficiently broad; for the particular
parameters used here the $0.025$ and $0.975$ centiles of the effective
prior distribution on $J(x;Q_y)$ are respectively -2 and +2 nepers
approximately. This might not be sufficiently wide for prediction
algorithms that provide much more information, and suitable
adjustments in the parameters in section \ref{infoparameters} will
then be needed, particularly if the number of data points is low.

Critics of this approach will surely point out that the results
therefore depend on the prior used for this modelling and the precise
samples $\theta_k$ used, both of which are true. More importantly,
there are likely to be some instances where the distribution of the
$j(x,y)$ carries a known parametric form (though we have not yet
encountered such a situation), in which circumstances it would clearly
be more appropriate to use that parametric form rather than one
introduced \textit{ex vacuo}.

\section{Discussion}
\label{conclusions}

In the situation where we are interested in predicting the value of a
variable $x$, with the prediction potentially to be used by many users
with different cost functions, the key attribute of a high quality
prediction algorithm is that its predictions should not be wrong,
i.e. that the true value of $x$ should be predicted with as high a
probability (or density) as possible.

While it is true that point value predictions may suffice where all
that matters is, for example, the mean squared distance of the point
prediction from the true value, there are many situations where the
relationship between the true value and its eventual effects is much
more complicated. Measurement of apparent Shannon information tells us
exactly how reliable the original prediction is compared with the
prior, and simultaneously tells us how reliable derived values are, at
least for smoothly and bijectively related ones.

As indicated in section \ref{proplist} points \ref{unique1} and
\ref{unique2}, ASI is characterised by being the unique quality
measure with a simple list of desirable properties; indeed there are
two such lists, of very different nature, each of which yields the
same quality measure.

The effect of optimising a prediction according to the ASI criterion
is often counter-intuitive; tests with naive human subjects usually
result in negative values even with respect to very broad priors, even
when attempting to predict values of things you might think the
subjects would know a lot about. It is also much more difficult to
obtain positive ASI from a prediction algorithm than one might think;
the characterisation by the betting scheme of section \ref{betting}
perhaps gives some intuitive idea why this might be so. This is why we
believe this is the right way to assess high quality prediction
algorithms -- it is sensitive to their potential failings.

A genuine argument against using the ASI in the past is that it has
been difficult to know how uncertain the value is that is obtained by
simply averaging values $j(x,y)$ over a finite previously unseen
dataset. The difficulty of doing this appears to have lain in the
asymmetry of the heavy tails of typical distributions of $j(x,y)$. We
believe that the remedy to this lies as indicated in extending the
Bayesian nature of the investigation process to cover not only the
prediction algorithm itself, but also the modelling of the uncertainty
of $j(x,y)$, and in particular in using the skew-Student distribution
as the fundamental building block of such a model. Clearly such a
process is affected by the prior used, but in practice, as in many
applications of Bayesian inference, so long as the prior is chosen to
be sufficiently broad and ``uninformative'', the exact details of the
prior turn out to be relatively unimportant. In this particular
instance it appears in particular to be important that the
``uninformative'' nature of the prior also covers the extent of
asymmetry, if any, in individual mixture component distributions.

\appendix

\section{Outline proof of uniqueness}
\label{appxproof}

In the following discussion we will assume for simplicity that all the
random variables are real-valued. However this is not a necessary
assumption, and the reader should have no difficulty converting for
example to the situation where random variables have values in
$\mathbb{R}^n$, where the main change needed is that $g'(x)$ (the
absolute value of the derivative of $g$ at $x$) needs replacing by
$\left|\det\left(\frac{\partial g}{\partial x}\right)\right|$, the
absolute value of the Jacobian determinant.

Suppose then $\phi(z_1,z_0,E)$ is a real-valued quality measure taking
as arguments two random variables and a conditioning event $E$ which
is unobserved but restricts us to only part of the total probability
space. Suppose $q_1,q_0$ are two functions of a real variable giving
density functions on the reals, so that if $x,y$ are random variables
then $q_i(y)(x)$ is also a random variable.

Intuitively, think of $q_1(y),q_0(y)$ as being the predictions of $x$
made by the two different algorithms from data $y$, and of
$\phi(q_1(y)(x),q_0(y)(x),E)$ as (somehow) measuring the quality of
algorithm represented by $q_1$ relative to that represented by $q_0$
under the circumstance that $E$ is true (but not known to be), so that
$q_1(r)$ denotes the probability distribution on $x$ resulting from
the algorithm trying to deduce the value of the random variable $x$
from that of $y$ given the particular value $r$ of $y$.

More technically, suppose $(\Omega,M,P)$ is a probability space, $X$
is the set of real-valued random variables on it having density
functions, $G$ is the set of diffeomorphisms\footnote{A diffeomorphism
here is a bijective differentiable function whose inverse is
differentiable.} $:\mathbb{R}\to \mathbb{R}$, $H$ is the set of
bijections $:\mathbb{R}\to \mathbb{R}$, $F$ is the set of probability
density functions on $\mathbb{R}$, and $Q$ the set of functions
$:\mathbb{R}\to F$.

Now suppose that $x,y\in X$ and that $\phi:X\times X \times M \to
\mathbb{R}$ is a function. In the following we use the notation
$P(x@a)$ to denote the value of the probability density function of
the random variable $x$ at the value $a$.

Suppose $\phi$ has the following properties for all $z_1,z_0\in
X,q,q_0,q_1,q_2\in Q, E \in M, g \in G$:

\begin{enumerate}

\item \label{prop1} If $\hat{q}(r)(s)$ is defined to be $P(x@s|y=r)$ and similarly
  $\bar{q}(r)(s)=P(x@s)$, then
  $$I(x;y)=\phi(\hat{q}(y)(x),\bar{q}(y)(x),\Omega).$$ (The value on
  the Bayesian posterior relative to the prior is the true Shannon
  mutual information.)

\item \label{prop2} For any $g\in G$, $$\phi(T_g(q_1)(y)(g^{-1}(x)),
  T_g(q_0)(y)(g^{-1}(x)),E) = \phi(q_1(y)(x),q_0(y)(x),E),$$ where
  $T_g(q)(r)(s)=g'(s)q(r)(g(s))$, with the dash denoting
  differentiation. (Consequent predictions of diffeomorphically
  related variables have the same prediction quality; $T_g$ transforms
  a prediction $q$ about $s$ given $r$ into one about $g^{-1}(s)$.)

\item \label{prop3} $$\phi(z_2,z_0,E) = \phi(z_2,z_1,E)+\phi(z_1,z_0,E).$$ (The
  increase in quality from algorithm 0 to algorithm 2 is the sum of
  the increases from 0 to 1 and from 1 to 2.)

\item \label{prop4} If $E_1,E_2\in M,E_1\cap E_2=\emptyset$ then
$$P(E_1\cup E_2) \phi(z_1,z_0,E_1\cup E_2) =
  P(E_1)\phi(z_1,z_0,E_1)+P(E_2)\phi(z_1,z_0,E_2).$$ (The quality
  measure averages correctly over unobserved other circumstances.)

\end{enumerate}

We then claim that it follows (assuming the existence of sufficiently
varied random variables) that $$\phi(q_1(y)(x),q_0(y)(x),E) =
\int{P(x,y|E)\log\left(\frac{q_1(y)(x)}{q_0(y)(x)}\right)\,d(x,y)},$$
and hence by \ref{prop3} any absolute quality measure that satisfies
these properties and assigns zero to the algorithm that just outputs
the prior distribution on $x$ and ignores the value of $y$ is the
apparent Shannon information.

We will however only sketch the outline of a proof, leaving the reader
to fill in the details; this outline goes as follows.

Property \ref{prop4} is first used to show that $\phi(z_1,z_0,E)$ is
the expectation given $E$ of a random variable $:\Omega\to\mathbb{R}$
which sends $\omega$ to a function of $(z_1(\omega),z_0(\omega))$, so
that we can write $$\phi(z_1,z_0,E) =
\int{P(z_1,z_0|E)f(z_1,z_0)\,d(z_1,z_0)}$$ for some function
$f:\mathbb{R}^2\to\mathbb{R}.$

Next, property \ref{prop2} allows us to deduce that for any positive
real $r$, $\phi(rz_1,rz_0,E) = \phi(z_1,z_0,E)$, so in fact
$f(z_1,z_0)$ is a function of $\alpha = \frac{z_1}{z_0}$. Rewriting
$f(z_1,z_0)$ as $f(\alpha)$, property \ref{prop3} then tells us that
$f(\alpha_1 \alpha_2) = f(\alpha_1)+f(\alpha_2)$, from which we learn
that $f$ is a constant multiple of the logarithm.

Finally, property \ref{prop1} tells us that the scalar multiple must
be unity, which gives the desired result.

\section{Details of modelling method used}
\label{modelling}

This appendix contains the details of the modelling used in section
\ref{approach} above. The case of modelling the $j(x,y)$ data was done
with one set of hyperparameter values (listed in section
\ref{infoparameters}), the case of modelling the recurrence times for
prostate cancer patients with a different set (listed in section
\ref{prostateparameters}), as well as some further modifications to
allow for some recurrence times being unobserved, which are listed in
section \ref{prostatemods} below.

\subsection{Overall model hierarchy}
\label{hierarchy}

Let $x_k$ denote the $k$th data point to be modelled.

Let $C$ denote the number of mixture components, $c$ an individual
component, $c_k$ the mixture component to which $x_k$ is adherent, and
$\mathbf{c}$ the vector of the $c_k$s.

For other variables subscripts are used for two purposes: first, to
indicate which component of a vector of such parameters is being used;
second, to modify the variable name to indicate a similar variable at
a higher level of the model. Thus $\mathbf{S}$ denotes the vector of
scale matrices for the mixture components, $\mathbf{S}_c$ denotes the
scale matrix for mixture component $c$, $\mathbf{R}_\mathbf{S}$
denotes the vector of scale matrices for $\mathbf{S}$, while
$\mathbf{R}_{\mathbf{S}_c}$ denotes the component of
$\mathbf{R}_\mathbf{S}$ that applies to $\mathbf{S}_c$

All variable with names starting $\mathbf{S}$ or $\mathbf{R}$ are
positive definite symmetric $N \times N$ real matrices.

Then the overall model hierarchy can be represented by the following
equations; we start with those at the top of the hierarchy.
$\mu_{\mu_\mu}, \mathbf{S}_{\mu_\mu}, m_{\mathbf{S}_\mu},
m_{\mathbf{R}_{\mathbf{S}_\mu}},
\mathbf{R}_{\mathbf{R}_{\mathbf{S}_\mu}}, N, \kappa_\nu, \kappa_\eta,
\kappa_C, a_m, b_m, N_m, m_{\mathbf{R}_\mathbf{S}},
\mathbf{R}_{\mathbf{R}_\mathbf{S}}, a_{m_\mathbf{S}},
b_{m_\mathbf{S}}$ are constant hyperparameters, values given in
sections \ref{infoparameters} or \ref{prostateparameters} below
depending on application.

$$P(C|\kappa_C) = \left\{\begin{matrix}\kappa_C^{C-1}(1-\kappa_C) &
    (C>0)\\
0 & (C=0)\end{matrix}\right.\ \ \ \text{(exponential on integers)}$$

$$\eta_c=\kappa_\eta / C$$

$$P(p|\eta,C) = \left\{\begin{matrix}\frac{1}{\sqrt{C}}
    \frac{\Gamma(\sum_{c=1}^C{\eta_c})}{\prod_{c=1}^C{\Gamma(\eta_c)}}
    \prod_{c=1}^C{p_c^{\eta_c-1}} &
      \left(\sum_{c=1}^C{p_c}=1\right)\\ 0 &
      \text{(otherwise)}\end{matrix}\right.\ \ \ \text{(Dirichlet)}$$

$$P(c_k@c|p)=p_c$$

$$P(\mathbf{R}_{\mathbf{S}_\mu} | m_{\mathbf{R}_{\mathbf{S}_\mu}}, \mathbf{R}_{\mathbf{R}_{\mathbf{S}_\mu}}, N)
 = \text{Wishart}(\mathbf{R}_{\mathbf{S}_\mu}; m_{\mathbf{R}_{\mathbf{S}_\mu}}, \mathbf{R}_{\mathbf{R}_{\mathbf{S}_\mu}}, N)$$
where 
$$\text{Wishart}(\mathbf{S};m,\mathbf{R},N) =
\frac{\det(\mathbf{R})^{m+\frac{N-1}{2}}}{(2\pi)^{N(N-1)/4}\prod_{j=0}^{N-1}{\Gamma\left(m+\frac{j}{2}\right)}}
\det(\mathbf{S})^{m-1} e^{-\trace(\mathbf{R}\mathbf{S})}$$

$$P(\mathbf{S}_\mu|m_{\mathbf{S}_\mu}, \mathbf{R}_{\mathbf{S}_\mu}, N)
= \text{Wishart}(\mathbf{S}_\mu;m_{\mathbf{S}_\mu}, \mathbf{R}_{\mathbf{S}_\mu}, N)$$

$$P(\mu_\mu|\mu_{\mu_\mu}, \mathbf{S}_{\mu_\mu}) =
\text{Gaussian}(\mu_\mu; \mu_{\mu_\mu}, \mathbf{S}_{\mu_\mu})$$
where
$$\text{Gaussian}(x;\mu,\mathbf{S}) =
\sqrt{\det\left(\frac{\mathbf{S}}{2\pi}\right)}
e^{-\frac{1}{2}(x-\mu)'\mathbf{S}(x-\mu)}$$
where $'$ denotes transpose.

$$P(m_\mathbf{S}|a_{m_\mathbf{S}}, b_{m_\mathbf{S}}, N) \propto
\text{proGamma}(m_\mathbf{S};a_{m_\mathbf{S}}, b_{m_\mathbf{S}}, N,
1)$$ where the proGamma distribution is defined below in section
\ref{proGamma}.

$$P(\mathbf{R}_\mathbf{S}|m_{\mathbf{R}_\mathbf{S}},
\mathbf{R}_{\mathbf{R}_\mathbf{S}}, N) = \text{Wishart}(\mathbf{R}_\mathbf{S};m_{\mathbf{R}_\mathbf{S}},
\mathbf{R}_{\mathbf{R}_\mathbf{S}}, N)$$

$$P(\mathbf{S}_c|m_\mathbf{S}, \mathbf{R}_\mathbf{S}, N) =
\text{Wishart}(\mathbf{S}_c;m_\mathbf{S},
(m_\mathbf{S}-1)\mathbf{R}_\mathbf{S}, N)$$

$$P(m_c|a_m,b_m,N_m)\propto \text{proGamma}(m_c;a_m,b_m,N_m,1)$$

$$r_c=m_c-1$$

$$P(\alpha_k|m,r) =
\frac{r_{c_k}^{m_{c_k}}}{\Gamma(m_{c_k})}\alpha_k^{m_{c_k}-1}e^{-r_{c_k}\alpha_k}
\ \ (\alpha_k>0)\ \ \ \text{(Gamma)}$$

$$P(\nu_c|\kappa_\nu,\mathbf{S}_c) = \text{Gaussian}\left(\nu_c;
0,\frac{\mathbf{S}_c}{\kappa_\nu}\right)$$

$$P(\mu_c|\mu_\mu,\mathbf{S}_\mu) =
\text{Gaussian}(\mu_c;\mu_\mu,\mathbf{S}_\mu)$$

$$P(x_k|\mathbf{c},\mu,\nu,\alpha,\mathbf{S}) =
\text{Gaussian}\left(x_k;
\mu_{c_k}+\frac{\nu_{c_k}}{\sqrt{\alpha_k}},\alpha_k\mathbf{S}_{c_k}\right)$$

For convenience we also include the detailed definitions of the
perhaps less familiar skew-Student and proGamma distributions in
sections \ref{skewStudent} and \ref{proGamma} below.

\subsection{Modifications to the model for modelling prostate cancer recurrence
times}
\label{prostatemods}

For modelling prostate cancer recurrence times, the above model was
modified in the following ways.

First, the scalar $x_k$ representing a $j(x,y)$ value is replaced by a
vector $x_k$ whose first component $x_{k,1}$ represents the logarithm
of the true time at which recurrence occurs (which may or may not be
observed), and whose remaining components represent logarithms of the
various biochemical marker levels and clinical variables used for the
prediction. 

Second, the observed data vector $y_k$ is modelled as follows:
$$y_{k,1} = (x_{k,1}+n_{k,1})\land z_k$$
$$y_{k,\bar{1}} = (x_{k,\bar{1}}+n_{k,\bar{1}})$$
where $a\land b$ denotes the minimum of $a$ and $b$ and $\bar{1}$
denotes all indices other than 1, and where 
$$P(n_k|m_n,r_n,\mathbf{S}_n,N) =
\sqrt{\det\left(\frac{\mathbf{S}_n}{2\pi}\right)}
\frac{\Gamma\left(m_n+\frac{N}{2}\right)}{\Gamma(m_n)}
\frac{r_n^{m_n}} {\left(r_n+\frac{1}{2}n_k'\mathbf{S}_n
  n_k\right)^{m_n+\frac{N}{2}}}\ \ \ \text{(Student)}$$

$$r_n=m_n$$

$$P(z_k)=\text{Gaussian}(z_k;\mu_z,\mathbf{S}_z)$$ where
$m_n,\mathbf{S}_n,\mu_z,\mathbf{S}_z$ are additional constant
hyperparameters, and $z_k$ is the logarithm of the duration of
observation if recurrence was not observed in case $k$ (and is known
to be greater than $y_{k,1}$ otherwise). $y_{k,1}$ denotes the earlier
of the time the patient's follow-up ceased with the patient
recurrence-free or the time cancer returned.

Third, the following equations replace the corresponding ones in the hierarchy of
section \ref{hierarchy}:
$$P(m_c|a_m,b_m,N_m) \propto \text{proGamma}(m_c;a_m,b_m,N_m,2)
\ \ \text{(proGamma type 2 – see section \ref{proGamma}))}$$
$$r_c=m_c$$

\subsection{Sampling from the posterior distribution of the model}

Exploration of the posterior distribution given the data was done
using Markov chain Monte-Carlo by Gibbs resampling along sets of
mutually independent axes, using standard methods. In particular the
proGamma distribution can be resampled using adaptive rejection
sampling\cite{GilksWild}.

Convergence time was assessed by creating synthetic datasets whose
true parameters were known and ensuring that run time was long enough
for the resulting distributions to be approximately independent of
starting parameters, where the starting parameters were initialised
either from the true values or from random values from their priors.
For modelling the distribution of $J(x;Q_y)$, having obtained many
sets of posterior joint parameter samples given the information values
$j(x,y)$ whose distribution is being modelled, the mean of each such
distribution was calculated analytically and considered as a set of
posterior samples of the apparent information content being
estimated. Mean, median, and centiles of this distribution were then
estimated from the respective statistics of this set of samples.

For modelling recurrence times of prostate cancer, having obtained
many sets of posterior joint parameter samples given the observed data
vectors, the marginal distribution of recurrence time was calculated
for each such joint parameter sample, and the equiprobable mixture of
these distributions returned as the predictive distribution on
recurrence time.

\subsection{Methods for choosing hyperparameters}

Bayesian model choice was undertaken for the establishment of
appropriate values of the top level hyperparameters, using
thermodynamic integration\cite{NealMCMC}. Where uncertainty of the
posterior model probabilities was too big to choose between close
candidate values of the hyperparameters, samples from the overall
prior model were examined by the human eye to look for similarity to
the observed data patterns.

\subsection{The skew-Student distribution}
\label{skewStudent}

The combined effect of the equations in section \ref{hierarchy} for
the distributions of $\alpha_k$ and $x_k$ is that \\$
P(x_k|\mathbf{c},\mu,\nu,m,r,\mathbf{S},N) =
\frac{r^m}{\Gamma(m)}\sqrt{\det\left(\frac{\mathbf{S}}{2\pi}\right)}
\frac{1}{\left(r+\frac{1}{2}\mathbf{y}'\mathbf{S}\mathbf{y}\right)^{m+\frac{N}{2}}}
e^{-\frac{1}{2}\nu'\mathbf{S}\nu}
\times\\ \left(\Gamma\left(m+\frac{N}{2}\right){_1F_1}\left(m+\frac{N}{2};\frac{1}{2};
\frac{(\mathbf{y}'\mathbf{S}\nu)^2}{4\left(r+\frac{1}{2}\mathbf{y'Sy}\right)}\right)
+\frac{\mathbf{y'S}\nu}{\sqrt{r+\frac{1}{2}\mathbf{y'Sy}}}\Gamma\left(m+\frac{N+1}{2}\right)
{_1F_1}\left(m+\frac{N+1}{2};\frac{3}{2};\frac{(\mathbf{y'S}\nu)^2}{4\left(r+\frac{1}{2}\mathbf{y'Sy}\right)}\right)\right)$
where $_1F_1$ denotes the relevant hypergeometric function, given
by $$_1F_1(a;b;z) =
\sum_{j=0}^\infty{\frac{\Gamma(a+j)/\Gamma(a)}{\Gamma(b+j)/\Gamma(b)}\frac{z^j}{j!}},$$
and where $\mathbf{y},\mu,\nu,m,\mathbf{R},\mathbf{S}$ respectively
stand for
$(\mathbf{x}_k-\mu_{c_k}),\mu_{c_k},\nu_{c_k},m_{c_k},\mathbf{R}_{c_k},\mathbf{S}_{c_k}$. The properties of the 1-dimensional
version of this distribution were reviewed in \cite{Shuford}.

\subsection{The proGamma distribution}
\label{proGamma}

(We would be glad to hear from anybody with a better known name for
this family of distributions, as also from anybody able to evaluate
the normalisation constants.)

The proGamma distributions are the conjugate distributions for the
shape parameter $m$ of the Gamma and Wishart distributions, with
respect to three different parameterisations of the Gamma
distribution.

The type 0 proGamma distribution is the simplest, conjugate to the
usual parameterisation of the Gamma with shape parameter $m$ and scale
parameter $r$, i.e. $$P(x|m,r)=\frac{r^m}{\Gamma(m)}x^{m-1}e^{-rx},$$
  as well as to the Wishart parameterised by $m$ and $\mathbf{R}$,
  i.e.
$$P(\mathbf{S}|m,\mathbf{R},N) =
\frac{\det(\mathbf{R})^{m+\frac{N-1}{2}}}{(2\pi)^{N(N-1)/4}\prod_{j=0}^{N-1}{\Gamma\left(m+\frac{j}{2}\right)}}
\det(\mathbf{S})^{m-1} e^{-\trace(\mathbf{R}\mathbf{S})}.$$ It is
given by 
$$P(m|a,b,N,v=0)\propto
\frac{e^{-(a+Nb)m}}{\left(\prod_{j=0}^{N-1}{\Gamma\left(m+\frac{j}{2}\right)}\right)^b}$$
  for $m>0$, where $N$ is 1 for the Gamma and is otherwise the
  dimensionality of the Wishart. Parameter restrictions are that $b$
  must be positive and $N$ must be a positive integer.

The type 1 proGamma distribution is instead conjugate with respect to
$n$ to the Gamma distribution parameterised by $(m,r) = (n,(n-1)t)$
(or for the Wishart $(m,\mathbf{R}) = (n,(n-1)\mathbf{T}))$, so that
the mean of the inverse of the Gamma or Wishart distributed quantity
doesn’t vary as $n$ varies. It now takes the form
$$P(m|a,b,N,v=1)\propto\left\{\begin{matrix}\frac{e^{-(a+Nb)m}(m-1)^{Nb\left(m+\frac{N-1}{2}\right)}}
{\left(\prod_{j=0}^{N-1}{\Gamma\left(m+\frac{j}{2}\right)}\right)^b} &
(m>1)\\ 0 & \text{(otherwise)}\end{matrix}\right.,$$ where the
parameter restrictions are again that both $a$ and $b$ must be
positive, but the support of the distribution is now $m > 1$ rather
than the $m > 0$ which applies for types 0 and 2.

The type 2 proGamma distribution is instead conjugate with respect to
$n$ to the Gamma parameterised by $(m,r) = (n,nt)$, so that as $n$
varies the mean of the Gamma remains constant. Similarly for the
Wishart it is conjugate with respect to $n$ where $(m,\mathbf{R}) =
\left(n, \left(n+\frac{N-1}{2}\right)\mathbf{T}\right),$ so that again
the mean of the Wishart remains constant as $n$ varies. It takes the
form $$P(m|a,b,N,v=2)\propto
\frac{e^{-(a+Nb)m}\left(m+\frac{N-1}{2}\right)^{Nb\left(m+\frac{N-1}{2}\right)}}
     {\left(\prod_{j=0}^{N-1}{\Gamma\left(m+\frac{j}{2}\right)}\right)^b}$$
     for $m>0$, where now the parameter restrictions are that both $a$
     and $b$ must be positive.

\subsection{Parameter settings for modelling of $j(x,y)$ data}
\label{infoparameters}

The values of the hyperparameters for the model of section
\ref{hierarchy} when used for modelling the uncertainty of $J(x;Q_y)$
based on samples of $j(x,y)$ are as follows:

\begin{center}
\begin{tabular}{cccc}
$N=1$ & $N_m=1$ & $\kappa_\nu=1$ & $\kappa_\eta=10$\\
$\kappa_C=0.9$ & $\mu_{\mu_\mu}=0$ & $\mathbf{S}_{\mu_\mu}=1$ &
  $m_{\mathbf{S}_\mu}=1.1$\\
$m_{\mathbf{R}_{\mathbf{S}_\mu}}=2$ &
  $\mathbf{R}_{\mathbf{R}_{\mathbf{S}_\mu}}=2.8$ & $a_m=1$ & $b_m=3$\\
$m_{\mathbf{R}_\mathbf{S}}=2$ &
  $\mathbf{R}_{\mathbf{R}_\mathbf{S}}=200$ & $a_{m_\mathbf{S}}=1$ &
  $b_{m_\mathbf{S}}=2$ 
\end{tabular}
\end{center}

\subsection{Parameter settings and model modifications for modelling of
prostate cancer recurrence times}
\label{prostateparameters}

The values of the hyperparameters for the model of section
\ref{hierarchy} as modified by section \ref{prostatemods} for
modelling the recurrence time of prostate cancer were as follows:

\begin{center}
\begin{tabular}{cccc}
$N=14$ & $N_m=1$ & $\kappa_\nu=0.1$ & $\kappa_\eta=10$\\ 
$\kappa_C=0.09$ & $\mu_{\mu_\mu}=\mean(T)$ &
  $\mathbf{S}_{\mu_\mu}=6.25\times(\cov(T))^{-1}$ &
  $m_{\mathbf{S}_\mu}=2.5$\\
$m_{\mathbf{R}_{\mathbf{S}_\mu}} = 2$ &
  $\mathbf{R}_{\mathbf{R}_{\mathbf{S}_\mu}} = 24 \times(\cov(T))^{-1}$
  & $a_m=1$ & $b_m=20$\\
$m_{\mathbf{R}_\mathbf{S}}=133$ &
  $\mathbf{R}_{\mathbf{R}_\mathbf{S}}=140\times(\cov(T))^{-1}$ &
  $a_{m_\mathbf{S}}=1$ & $b_{m_\mathbf{S}}=0.01$\\
$m_n=5$ & $\mathbf{S}_n=100 I_N$ & $\mu_z=\mu_{\mu_\mu,1}$ &
  $\mathbf{S}_z=(\cov(T)_{1,1})^{-1}$ 
\end{tabular}
\end{center}
where $T$ denotes the training dataset.

\bibliography{ms}
\bibliographystyle{ieeetr}

\end{document}